\documentclass[10pt]{article}
\usepackage{amsmath,amssymb,amsfonts}

\newtheorem{theorem}{Theorem}
\newtheorem{lemma}{Lemma}

\date{}

\title{Sampling Parts of Random Integer Partitions: A Probabilistic and Asymptotic Analysis}
\bigskip

\author{{\bf Ljuben Mutafchiev}\\
American University in Bulgaria, 2700 Blagoevgrad, Bulgaria \\ and
Institute of Mathematics and Informatics of the \\ Bulgarian
Academy of Sciences
\\ \tt {e-mail: ljuben@aubg.bg; tel.: +359 73 888498}}

\begin{document}
\maketitle

\begin{abstract}
Let $\lambda$ be a partition of the positive integer $n$, selected
 uniformly at random among all such partitions. Corteel et al.
 (1999) proposed three different procedures of sampling parts of
 $\lambda$ at random. They obtained limiting distributions of the
 multiplicity $\mu_n=\mu_n(\lambda)$ of the randomly-chosen part as $n\to\infty$.
 The asymptotic behavior of the part size
 $\sigma_n=\sigma_n(\lambda)$, under these sampling conditions was
 found by Fristedt (1993) and Mutafchiev (2014). All these results
 motivated us to study the relationship between the size and the
 multiplicity of a randomly-selected part of a random partition.
 We describe it obtaining the joint limiting distributions of
 $(\mu_n,\sigma_n)$, as $n\to\infty$, for all these three sampling procedures.
 It turns out that different sampling plans lead to different
 limiting distributions for $(\mu_n,\sigma_n)$. Our results
 generalize those obtained earlier and confirm the known
 expressions for the marginal limiting distributions of $\mu_n$
 and $\sigma_n$.
\end{abstract}

\vspace {.5cm}

\vspace{.5cm}


\vspace{.2cm}

\section{Introduction}

Partitioning integers into summands (parts) is a subject of
intensive research in combinatorics, number theory and statistical
physics. If $n$ is a positive integer, then by a partition,
$\lambda$, of $n$, we mean a representation
\begin{equation}\label{partition}
 \lambda: \quad n=\sum_{j=1}^n jm_j,
\end{equation}
 in which $m_j$, called multiplicities of parts $j,
j=1,2,...,n$, are non-negative integers. We use $\Lambda(n)$ to
denote the set of all partitions of $n$ and let
$p(n)=\mid\Lambda(n)\mid$. The number $p(n)$ is determined
asymptotically by the famous partition formula of Hardy and
Ramanujan [9]:
 \begin{equation}\label{hardy}
 p(n)\sim\frac{1}{4n\sqrt{3}}
\exp{\left(\pi\sqrt{\frac{2n}{3}}\right)}, \quad n\to\infty.
 \end{equation}
 A precise asymptotic expansion for $p(n)$ was found later by
Rademacher [14] (more details may be also found in [2; Chapter
5]). For instance, Rademacher's result implies that
\begin{eqnarray}\label{rademacher}
 & & p(n)=\frac{1}{4n\sqrt{3}}
\exp{\left(\pi\sqrt{\frac{2n}{3}}\right)}
-\frac{1}{4\pi\sqrt{2}n^{3/2}}\exp{\left(\pi\sqrt{\frac{2n}{3}}\right)}
\nonumber \\
 & & +O\left(\exp{\left(\frac{\pi}{2}\sqrt{\frac{2n}{3}}\right)}\right),
\quad n\to\infty.
\end{eqnarray}

 Further on, we assume that, for fixed integer $n\ge
1$, a partition $\lambda\in\Lambda(n)$ is selected uniformly at
random (uar), i.e. with probability $1/p(n)$. In this way, each
numerical characteristic of $\lambda$ can be regarded as a random
variable defined on the space $\Lambda(n)$.

Corteel et al. [3] proposed and studied three procedures of
sampling parts of a random partition $\lambda\in\Lambda(n)$. Basic
statistics of a randomly selected part are the part size and its
multiplicity. Corteel et al. [3] focused on the multiplicity
$\mu_{n,j}=\mu_{n,j}(\lambda)$ ($j=1,2,3,$) of the
randomly-selected part and found limiting distributions for
$\mu_{n,j}$, as $n\to\infty$, in these three cases of sampling
(here the subscript $j$ specifies the concrete sampling procedure
that is followed; the definitions of these three sampling
procedures will be given in the next section). In the same way,
let $\sigma_{n,j}=\sigma_{n,j}(\lambda)$ ($j=1,2,3$) be the size
of the randomly-selected part. Limit theorems for $\sigma_{n,j}$
were obtained in [6] and [13]. All these results motivated us to
study the relationship between the size and the multiplicity of a
randomly-selected part of a random integer partition. We describe
it obtaining the joint limiting distributions of $\mu_{n,j}$ and
$\sigma_{n,j}$ (j=1,2,3) as $n\to\infty$. Our results generalize
those obtained earlier in [6,3,13] and confirm the known
expressions for the marginal limiting distributions of $\mu_{n,j}$
and $\sigma_{n.j}$.

We organize our paper as follows. In Section 2 we describe the
sampling procedures proposed by Corteel et al. [3]. The main
results of this paper are stated in Section 3. The method of proof
is also briefly described there. Section 4 contains some auxiliary
facts on generating functions and some asymptotics that we need
further. We present the proofs of our limit theorems in Sections 5
- 7.

\section{Basic Random Variables and Definitions of the Sampling Procedures}
\setcounter{equation}{0}

  For any $\lambda\in\Lambda(n)$ selected uar, we define the random
 variables
 $$
 \alpha_j^{(n)} =\alpha_j^{(n)}(\lambda) =\mbox{the number of
 parts of size } j \mbox { in } \lambda.
 $$
 By $I_A$ we denote the indicator of an event $A$ and, for any two
 real numbers $d,s\ge 1$ and integer $m\ge 1$, we set
 \begin{equation}\label{zds}
 Z_{d,s}^{(n)}=\sum_{1\le j\le s}\alpha_j^{(n)}
 I_{\{\alpha_j^{(n)}\le d\}},
 \end{equation}
 \begin{equation}\label{yms}
 Y_{m,s}^{(n)}=\sum_{1\le j\le s} I_{\{\alpha_j^{(n)}=m\}}.
 \end{equation}
 ($Z_{d,s}^{(n)}$ counts the number of parts of size not grater than $s$ and multiplicity
 not greater than $d$ in a randomly-chosen partition $\lambda$, while $Y_{m,s}^{(n)}$ is the number of distinct parts
 with multiplicity $m$ and size not greater than $s$). Obviously,
 \begin{equation}\label{zn}
 Z_n=\sum_{j=1}^n \alpha_j^{(n)}
 \end{equation}
  equals the
 total number of parts and
 \begin{equation}\label{yn}
 Y_n =\sum_{j=1}^n I_{\{\alpha_j^{(n)}>0\}}
 \end{equation}
  - the
 number of distinct parts in $\lambda\in\Lambda(n)$.

 To describe the sampling procedures introduced by Corteel et al. [3] we
notice that they are two-step procedures that combine the outcomes
of two experiments. Therefore, they lead to three different
product  probability spaces. Since in each procedure we first
sample uar a partition $\lambda\in\Lambda(n)$, the probability
space on $\Lambda(n)$, equipped with the uniform probability
measure $Pr(\lambda\in\Lambda(n))=1/p(n)$, is included in each
product space. The second steps of sampling are, however,
different and therefore, for each different procedure we obtain a
different product space and different product probability measure.
In what follows next, we adopt the common notation $\mathbb{P}(.)$
for the product probability measure of each sampling procedure and
follow the concept of a product space developed in [8; Chapter
1.6]. By $\mathbb{E}(X)$ we denote the expected value of the
random variable $X$ defined on the integer partition space
$\Lambda(n)$.

{\it Procedure 1.} Given a partition $\lambda\in\Lambda(n)$ chosen
uar (step 1), we select a part uar among all $Z_n$ parts of
$\lambda$ (without any bias, step 2). By the product measure
formula [8; Chapter 1.6], (\ref{zds}) and (\ref{zn})
 \begin{eqnarray}
 & & \mathbb{P}(\{\lambda\in\Lambda(n)\}\times\{\mu_{n,1}\le d, \sigma_{n,1}\le
 s\}) \nonumber \\
 & & =Pr(\lambda\in\Lambda(n)) \mathbb{P}(\mu_{n,1}\le d, \sigma_{n,1}\le
 s)
  =\left(\frac{1}{p(n)}\right)\left(\frac{Z_{d,s}^{(n)}}{Z_n}\right). \nonumber
 \end{eqnarray}
Summation over all $\lambda\in\Lambda(n)$ yields
\begin{equation}\label{procone}
\mathbb{P}(\mu_{n,1}\le d, \sigma_{n,1}\le s)
=\mathbb{E}\left(\frac{Z_{d,s}^{(n)}}{Z_n}\right).
\end{equation}

{\it Procedure 2.} Given a partition $\lambda\in\Lambda(n)$ chosen
uar (step 1), we select a part among all $Y_n$ different parts
(step 2). Recalling definitions (\ref{yms}) and (\ref{yn}) of the
random variables $Y_{m,s}^{(n)}$ and $Y_n$, respectively, we
obtain in a similar way that
$$
\mathbb{P}(\{\lambda\in\Lambda(n)\}\times\{\mu_{n,2}=m,
\sigma_{n,2} \le s\})
=\left(\frac{1}{p(n)}\right)\left(\frac{Y_{m,s}^{(n)}}{Y_n}\right)
$$
and
\begin{equation}\label{proctwo}
\mathbb{P}(\mu_{n,2}=m, \sigma_{n,2}\le s)
=\mathbb{E}\left(\frac{Y_{m,s}^{(n)}}{Y_n}\right).
\end{equation}

{\it Procedure 3.} Given a partition $\lambda\in\Lambda(n)$ chosen
uar (step 1), we select a part of $\lambda$ with the probability
proportional to its size and multiplicity (step 2). Thus we set
\begin{equation} \label{thr}
\mathbb{P}(\{\lambda\in\Lambda(n)\}\times\{\mu_{n,3}=m,
\sigma_{n,3}\le s\})
 = \left(\frac{1}{p(n)}\right)\left(\frac{m}{n}\right) \sum_{1\le j\le
 s} jI_{\{\alpha_j^{(n)}=m\}},
\end{equation}
which in turn implies that
\begin{equation}\label{procthree}
 \mathbb{P}(\mu_{n,3}=m,
\sigma_{n,3}\le s) =\frac{m}{n}\sum_{1\le j\le s}
jPr(\alpha_j^{(n)}=m).
\end{equation}

{\it Remark.} Sampling procedure 3 can be interpreted in terms of
Ferrers diagrams - the graphical representations of the integer
partitions $\lambda\in\Lambda(n)$ [2; Chapter 1.3]. It is obtained
as follows. We use the notation $\lambda_k$ to denote the $k$th
largest part of $\lambda$ for $k$ a positive integer; if the
number of parts $Z_n$ of $\lambda$ is $<k$, then $\lambda_k=0$.
The Ferrers diagram illustrates (\ref{partition}) by a
two-dimensional array of dots, composed by $\lambda_1$ dots in the
first (most left) row, $\lambda_2$ dots in the second row, ...,
$\lambda_{Z_n}$ dots in the last $Z_n$th row. Therefore, a Ferrers
diagram may be considered as a union of disjoint blocks
(rectangles) of dots with base $j$ and height $\alpha_j^{(n)}$
(the multiplicity of part $j$). So, (\ref{thr}) and
(\ref{procthree}) imply that the sampling probability in Procedure
3 is proportional to the area of the block to which the chosen
part belongs.

\section{Statement of the Main Results and Brief Description of the Method of Proof}
 \setcounter{equation}{0}

For sampling procedures 1 - 3, we have proved the following limit
theorems.

\begin{theorem}
For the reals $u$ and $v$, we let
$$
 F(u,v) =\left\{\begin{array}{ll} 0 & \quad  \mbox {if} \quad \min{\{u,v\}}\le
 0 \\
 0 & \quad \mbox {if}\quad \min{\{u,v\}}>0 \quad \mbox {but}\quad u+v\le 1, \\
 u+v-1 & \quad \mbox {if}\quad 0<u\le 1, 0<v\le 1 \quad \mbox {and} \quad
 u+v>1, \\
 \min{\{1,v\}} & \quad \mbox {if} \quad u>1 \quad \mbox {and}
 \quad 0<v\le 1, \\
 \min{\{1,u\}} & \quad \mbox {if} \quad v>1 \quad \mbox {and}
 \quad 0<u\le 1, \\
 1 & \quad \mbox {if} \quad u>1 \quad \mbox {and} \quad v>1.
 \end{array}\right.
 $$
 Then, we have
 $$
 \lim_{n\to\infty}\mathbb{P}\left(\frac{2\log{\mu_{n,1}}}{\log{n}}\le
 u,\frac{2\log{\sigma_{n,1}}}{\log{n}}\le v\right) =F(u,v).
 $$
\end{theorem}

\begin{theorem} Let $0<t<\infty$. Then, for any positive integer
$m$, we have
$$
\lim_{n\to\infty}\mathbb{P}\left(\mu_{n,2}=m,
\frac{\pi\sigma_{n,2}}{\sqrt{6n}}\le t\right) =\int_0^t
e^{-my}(1-e^{-y})dy.
$$
\end{theorem}

\begin{theorem} Let $0<t<\infty$. Then, for any positive integer
$m$, we have
$$
\lim_{n\to\infty}\mathbb{P}\left(\mu_{n,3}=m,
\frac{\pi\sigma_{n,3}}{\sqrt{6n}}\le t\right)
=\frac{6m}{\pi^2}\int_0^t y(1-e^{-y})e^{-my}dy.
$$
\end{theorem}

{\it Remark 1.} Since the inequalities
$\frac{2\log{\mu_{n,1}}}{\log{n}}\le
 u,\frac{2\log{\sigma_{n,1}}}{\log{n}}\le v$ are equivalent to
 $\mu_{n,1}\le n^{u/2}, \sigma_{n,1}\le n^{v/2}$, respectively, Theorem 1
 implies that the proportion of parts of size $\le n^{v/2}$ and multiplicity
$\le n^{u/2}$, $0<u,v<1$, is approximately equal to $u+v-1$ if
$u+v>1$; if $u+v\le 1$ this proportion approaches zero as
$n\to\infty$. For the other two sampling procedures, Theorems 2
and 3 show that typically chosen part sizes are of order
$const\sqrt{n}$, while their multiplicities are finite - both
converge weakly to discrete random variables whose support is the
set $\{1,2,...\}$.

{\it Remark 2.} For the sake of completeness, we present here a
list of the known marginal limiting distributions for the size and
multiplicity of the randomly-chosen part. They can be obtained as
corollaries of Theorems 1-3. Proper references are also given.

$$
\lim_{n\to\infty}\mathbb{P}\left(\frac{2\log{\mu_{n,1}}}{\log{n}}\le
t\right)=t, \quad 0<t<1
$$
[3; p.195];

$$
\lim_{n\to\infty}\mathbb{P}\left(\frac{2\log{\sigma_{n,1}}}{\log{n}}\le
t\right)=t, \quad 0<t<1
$$
[6; p.712];

$$
\lim_{n\to\infty}\mathbb{P}(\mu_{n,2}=m)=\frac{1}{m(m+1)}, \quad
m=1,2,..
$$
[3; p.192];

$$
\lim_{n\to\infty}\mathbb{P}\left(\frac{\pi\sigma_{n,2}}{\sqrt{6n}}
\le t\right) =1-e^{-t} \quad 0<t<\infty
$$
[13; Theorem 2];

$$
\lim_{n\to\infty}\mathbb{P}(\mu_{n,2}=m) =\frac{6(2m+1)}{\pi^2
m(m+1)^2}, \quad m=1,2,...
$$
[3; p. 195];

$$
\lim_{n\to\infty}\mathbb{P}\left(\frac{\pi\sigma_{n,3}}{\sqrt{6n}}
\le t\right) =\frac{6}{\pi^2}\int_0^t\frac{y}{e^y-1}dy, \quad
0<t<\infty
$$
[13; Theorem 3].

  We conclude this section with a description of our method of proof.
  It combines probabilistic with analytical tools.
We employ Fristedt's conditioning device [6], which allows to
transfer probability distributions of linear combinations of the
multiplicities $\alpha_j^{(n)}$ into conditional distributions of
the corresponding linear combinations of independent and
geometrically distributed random variables. Using this method, we
show that, as $n\to\infty$, the expected values in (\ref{procone})
and (\ref{proctwo}) are close to the ratios of the expectations of
the random variables that are involved there. The asymptotic
behavior of the expectations of $Y_n$ and $Z_n$, defined by
(\ref{yn}) and (\ref{zn}), respectively, is well known:
\begin{equation}\label{eyn}
\mathbb{E}(Y_n)\sim\sqrt{6n}/\pi,
\end{equation}
\begin{equation}\label{ezn}
\mathbb{E}(Z_n)\sim(\sqrt{6n}/2\pi)\log{n}
\end{equation}
(see [16] and [3], respectively). We use combinatorial enumeration
identities for generating functions, Cauchy coefficient formula
and the saddle-point method in terms of Hayman admissibility
theory [10] (see also [5; Chapter VIII.5]) to obtain the
asymptotic behavior of $\mathbb{E}(Z_{d,s}^{(n)})$ (see
(\ref{procone})). Finally, (\ref{proctwo})) and (\ref{procthree})
are analyzed using an approach developed by Corteel et al. [3] and
based on Euler-MacLaurin sum formula.

\section{Generating Functions and the Analytical Background of the
Proofs}
 \setcounter{equation}{0}

We start with the notation $g(x)$ for the generating function of
 the sequence $\{p(n)\}_{n\ge 1}$. For $\mid x\mid<1$, $g(x)$
 admits the well known representation
 \begin{equation}\label{euler}
 g(x)=1+\sum_{n=1}^\infty p(n)x^n =\prod_{k=1}^\infty (1-x^k)^{-1}
 \end{equation}
 (see e.g. [2; Theorem 1.1]). Our first lemma is related to the
 probability generating function and the expectation of the random
 variable $Z_{d,s}^{(n)}$, defined by (\ref{zds}).

 \begin{lemma} For any reals $d, s\ge 1$ and complex variables $x$
 and $z$, satisfying $\mid x\mid<1$ and $\mid z\mid<1$, we have
 \begin{equation}\label{idpgfzds}
 1+\sum_{n=1}^\infty p(n)x^n \mathbb{E}(z^{Z_{d,s}^{(n)}})
 =g(x)\prod_{1\le j\le s} \frac{(1-(zx^j)^{d+1})(1-x^j)}{1-zx^j}.
 \end{equation}
 Moreover
 \begin{eqnarray}\label{idezds}
 & & 1+\sum_{n=1}^\infty p(n)x^n \mathbb{E}(Z_{d,s}^{(n)})
  \\
 & & =g(x)\left(\sum_{1\le j\le s}\frac{x^j}{1-x^j}
 -(d+1)\sum_{1\le j\le s} \frac{x^{j(d+1)}}{1-x^{j(d+1)}}\right)
 \prod_{1\le j\le s} (1-x^{j(d+1)}).\nonumber
 \end{eqnarray}
\end{lemma}

{\it Proof.} The generating function identity (\ref{idpgfzds})
follows from a more general argument developed in [15; Chapter
V.5]. To state it we need some preliminary notations. We let
$B\subset\{1,2,...\}$ and let $\Omega_j\subset\mathcal{N}_0
=\{0,1,...\}, j\ge 1,$ be a sequence of sets. By
$\widetilde{\sum}$ we denote a sum over all $j\in B$, satisfying
(\ref{partition}) with $m_j\in\Omega_j, j\ge 1$. Then, we have
\begin{equation}\label{sachkov}
\prod_{j\in B}\sum_{m_j\in\Omega_j} (z_j x^j)^{m_j} =1+\sum_{n\ge
1} x^n \widetilde{\sum} z_1^{m_1}z_2^{m_2}...z_n^{m_n},
\end{equation}
where $x, z_1, z_2,...$ are formal variables. In (\ref{sachkov})
we set $B=\{1,2,...,[s]\}$,
 $$
\Omega_j= \left\{\begin{array}{ll} \{0,1,...,[d]\} & \qquad  \mbox {if}\qquad j\le s, \\
 \mathcal{N}_0 & \qquad \mbox {if}\qquad j>s
 \end{array}\right.
 $$
 and
$$
x_j= \left\{\begin{array}{ll} x & \qquad  \mbox {if}\qquad j\le s, \\
 1 & \qquad \mbox {if}\qquad j>s.
 \end{array}\right.
 $$
 (Here $[s]$ and $[d]$ denote the integer parts of $s$ and $d$, respectively.)
 The required identity (\ref{idpgfzds}) now follows from
 (\ref{zds}) and (\ref{euler}). A differentiation with respect to
 $z$ in (\ref{idpgfzds}) leads to the the expectations of
 $Z_{d,s}^{(n)}$ and identity (\ref{idezds}).  $\rule{2mm}{2mm}$

 The next lemma establishes a similar generating function identity
 for the random variable $Y_{m,s}^{(n)}$ defined by (\ref{yms}).
 It can be proved repeating the argument from [3; Theorem 1].

 \begin{lemma} For any real number $s\ge 1$, positive integer
 $m$ and complex variables $x$ and $z$, satisfying $\mid x\mid<1$ and $\mid z\mid<1$, we have
 $$
 1+\sum_{n\ge 1, j\ge 0} x^n \mathbb{E}(z^{Y_{m,s}^{(n)}})
 =g(x)\prod_{1\le k\le s}(1+(z-1)x^{mk}(1-x^k)).
 $$
 This in turn implies that
 \begin{equation}\label{eyms}
 \mathbb{E}(Y_{m,s}^{(n)}) =\sum_{1\le k\le s}
 (p(n-mk)-p(n-(m+1)k)).
 \end{equation}
 \end{lemma}

 Further on, for the sake of simplicity, we let
 \begin{equation}\label{c}
 c=\frac{\pi}{\sqrt{6}}.
 \end{equation}
 We notice that Hardy-Ramanujan-Rademacher's formula in its form (\ref{rademacher})
 implies that
 $$
 p(n) =\frac{e^{2cn^{1/2}}}{4\sqrt{3}n}(1+O(n^{-1/2})), \quad
 n\to\infty.
 $$
 Using this expression, Corteel et al. [3; p. 190] have obtained
 the following asymptotic estimates.

 \begin{lemma} For enough large $n$, we have
 $$
 \frac{p(n-mj)}{p(n)} =\left(1+O\left(\frac{mj}{n^{3/2}}\right)
 +O((n-mj)^{-1/2})\right) e^{-cmj/n^{1/2}}
 $$
 \begin{equation}\label{ratio}
  = \left\{\begin{array}{ll} (1+O(n^{-1/2}))e^{-cmj/n^{1/2}} & \qquad  \mbox {if}\qquad mj\le n/2,  \\
 O(e^{-cn^{1/2}/2}) & \qquad \mbox {if}\qquad mj>n/2.
 \end{array}\right.
 \end{equation}
 \end{lemma}

 Lemma 3 enables us to interpret the sum in (\ref{eyms}) as
 a Riemann integral sum.

 Our next preliminary fact is related to Hardy-Ramanujan formula
 (\ref{hardy}). We shall present it into a slightly different
 form, which will be used further to find the asymptotic of $\mathbb{E}(Z_{d,s}^{(n)})$.
 To introduce the reader into the subject, we notice that
 Hardy-Ramanujan formula has been subsequently
 generalized in various directions most notably by Meinardus [11]
 (see also [2; Chapter 6]). Meinardus obtained
the asymptotic of the Taylor coefficients of infinite products of
the form
\begin{equation} \label{product}
\prod_{k=1}^\infty (1-x^k)^{-b_k}
\end{equation}
under certain general assumptions on the sequence of non-negative
numbers $\{b_k\}_{k\ge 1}$. Meinardus approach is based on
considering the Dirichlet generating series
\begin{equation} \label{diofzi}
D(z)=\sum_{k=1}^\infty b_k k^{-z}, \quad z=u+iv.
\end{equation}
Since we shall use this result, below we briefly describe
Meinardus assumptions avoiding their precise statements as well as
some extra notations and concepts. The first Meinardus assumption
($M_1$) specifies the domain $\mathcal{H}=\{z: \Re{(z)}= u\ge
-C_0\}, 0<C_0<1,$ in the complex plane, in which $D(z)$ has an
analytic continuation. The second one ($M_2$) is related to the
asymptotic behavior of $D(z)$, whenever $\mid\Im{(z)}\mid =\mid
v\mid\to\infty$. A function of the complex variable $z$ which is
bounded by $O(\mid\Im(z)\mid^{C_1}), 0<C_1<\infty$, in certain
domain of the complex plane is called function of finite order.
Meinardus second condition ($M_2$) requires that $D(z)$ is of
finite order in the whole domain $\mathcal{H}$. Finally, the
Meinardus third condition ($M_3$) implies a bound on the ordinary
generating function of the sequence $\{b_k\}_{k\ge 1}$. It can be
stated in a way simpler than the Meinardus original expression by
the inequality
$$
\sum_{k=1}^\infty b_k e^{-k\omega}\sin^2{(\pi ku)} \ge
C_2\omega^{-\epsilon_1}, \quad 0<\frac{\omega}{2\pi}<\mid
u\mid<\frac{1}{2},
$$
for sufficiently small $\omega$ and some constants $C_2,
\epsilon_1>0$ ($C_2=C_2(\epsilon_1)$) (see [7; p. 310]).

It is known that Euler partition generating function $g(x)$ (which
is obviously of the form (\ref{product})) satisfies the Meinardus
scheme of conditions ($M_1$)-($M_3$) (see e.g. [2; Theorem 6.3]).

The proof of our Theorem 1 will be based on an asymptotic analysis
of a Cauchy integral stemming from (\ref{idezds}). We shall apply
there the saddle-point method in the sense of Hayman [10] (see
also [5; Chapter VIII.5]). In [10] Hayman studied a wide class of
power series satisfying a set of relatively mild conditions and
established general formulas for the asymptotic order of their
coefficients. In the proof of Theorem 1 we shall essentially use
that the generating function $g(x)$ is admissible in the sense of
Hayman.  To present Hayman's idea and show how it can be applied,
we need to introduce some auxiliary notations.

We consider here a function $G(x)=\sum_{n=1}^\infty G_n x^n$ that
is analytic for $\mid x\mid<\rho, 0<\rho<\infty$. For $0<r<\rho$,
we let
 \begin{equation}\label{ar}
 a(r)=r\frac{G^\prime(r)}{G(r)},
\end{equation}
\begin{equation}\label{br}
b(r)= r\frac{G^\prime(r)}{G(r)}
+r^2\frac{G^{\prime\prime}(r)}{G(r)}
-r^2\left(\frac{G^\prime(r)}{G(r)}\right).
\end{equation}
In the statement of Hayman's result we use the terminology given
in [5; Chapter VIII.5]. We assume that $G(x)>0$ for $x\in
(R_0,\rho)\subset (0,\rho)$ and satisfies the following three
conditions.

{\it Capture condition.} $\lim_{r\to\rho} a(r)=\infty$ and
$\lim_{r\to\rho} b(r)=\infty$.

{\it Locality condition.} For some function $\delta=\delta(r)$
defined over $(R_0,\rho)$ and satisfying $0<\delta<\pi$, one has
$$
G(re^{i\theta})\sim G(r)e^{i\theta a(r)-\theta^2 b(r)/2}
$$
as $r\to\rho$, uniformly for $\mid\theta\mid\le\delta(r)$.

{\it Decay condition.}
$$
G(re^{i\theta}) =o\left(\frac{G(r)}{\sqrt{b(r)}}\right)
$$
as $r\to\rho$, uniformly for $\delta(r)\le\theta<\pi$.

{\textbf{Hayman Theorem.} Let $G(x)$ be Hayman admissible function
and $r=r_n$ be the unique solution in the interval $(R_0,\rho)$ of
the equation
\begin{equation}\label{eqar}
a(r)=n.
\end{equation}
Then the Taylor coefficients of $G(x)$ satisfy, as $n\to\infty$,
\begin{equation}\label{hayman}
G_n\sim\frac{G(r_n)}{r_n^n\sqrt{2\pi b(r_n)}}
\end{equation}
 with $b(r_n)$ given
by (\ref{br}).

The next lemma presents an alternative formula for the partition
function $p(n)$.

\begin{lemma} If $r=r_n$ satisfies (\ref{eqar}) for
sufficiently large $n$, then
$$
p(n)\sim\frac{r_n^n g(r_n)}{\sqrt{2\pi b(r_n)}}, \quad n\to\infty,
$$
where $a(r_n)$ and $b(r_n)$ are given by (\ref{ar})and (\ref{br})
with $G(x)\equiv g(x)$.
\end{lemma}

{\it Proof.} Since in (\ref{euler}) we have $b_k=1, k\ge 1$, the
Dirichlet generating series (\ref{diofzi}) is $D(z)=\zeta(z)$,
where $\zeta$ denotes the Riemann zeta function. We set in
(\ref{ar}) and (\ref{br}) $r=r_n=e^{-h_n}, h_n>0$, where $h_n$ is
the unique solution of the equation
\begin{equation}\label{maineq}
a(e^{-h_n})=n.
\end{equation}
((\ref{maineq}) is an obvious modification of (\ref{eqar}).)
Granovsky et al. [7] showed that the first two Meinardus
conditions imply that the unique solution of (\ref{maineq}) has
the following asymptotic expansion:
\begin{equation}\label{dn}
h_n =\sqrt{\zeta(2)/n} +\frac{\zeta(0)}{2n} +O(n^{-1-\beta})
=\frac{\pi}{\sqrt{6n}}-\frac{1}{4n} +O(n^{-1-\beta}),
\end{equation}
where $\beta>0$ is fixed constant (here we have also used that
$\zeta(0)=-1/2$; see [1; Chapter 23.2]). We also notice that
(\ref{br}) and (\ref{dn}) impliy that
\begin{equation}\label{bdn}
b(e^{-h_n}) =2\zeta(2)h_n^{-3} +O(h_n^{-2}) \sim\frac{\pi^2}{3}
h_n^{-3} \sim\frac{2\sqrt{6}}{\pi} n^{3/2}
\end{equation}
(see [12; Lemma 2.2] with $D(z)=\zeta(z)$). Hence, by
(\ref{maineq}) and (\ref{bdn}), $a(e^{-h_n})\to\infty$ and
$b(e^{-h_n})\to\infty$ as $n\to\infty$, that is, Hayman's
``capture'' condition is satisfied with $r=r_n=e^{-h_n}$. To show
next that Hayman's ``decay'' condition is satisfied by $g(x)$ we
set
\begin{equation}\label{delta}
\delta_n =\frac{h_n^{4/3}}{\Omega(n)}
=\frac{\pi^{4/3}}{(6n)^{2/3}\Omega(n)}
\left(1+O\left(\frac{1}{\sqrt{n}}\right)\right)
\end{equation}
with $h_n$ given by (\ref{dn}), where $\Omega(n)\to\infty$ as
$n\to\infty$ arbitrarily slowly. We can apply now an estimate for
$\mid g(e^{-h_n+i\theta})\mid$ established in a general form in
[12; Lemma 2.4] using all three Meinardus conditions. It states
that there are two positive constants $c_0$ and $\epsilon_0$, such
that, for sufficiently large $n$,
\begin{equation}\label{decay}
\mid g(e^{-h_n+i\theta})\mid \le
g(e^{-h_n})e^{-c_0h_n^{-\epsilon_0}}
\end{equation}
uniformly for $\delta_n\le\mid\theta\mid<\pi$. This, in
combination with (\ref{bdn}), implies that $\mid
g(e^{-h_n+i\theta})\mid=o(g(e^{-h_n})/\sqrt{b(e^{-h_n})})$
uniformly in the same range for $\theta$, which is just Hayman's
``decay'' condition. Finally, by Lemma 2.3 of [12], established
using Meinardus conditions ($M_1$) and ($M_2$), Hayman's
``locality'' condition is also satisfied by $g(x)$. In fact, this
lemma implies in the particular case $D(z)=\zeta(z)$ that
\begin{equation}\label{locality}
e^{-i\theta n}\frac{g(e^{-h_n+i\theta})}{g(e^{-h_n})}
=e^{-\theta^2 b(e^{-h_n})/2}(1+O(1/\Omega^3(n))
\end{equation}
uniformly for $\mid\theta\mid\le\delta_n$, where $b(e^{-h_n})$ and
$\delta_n$ are determined by (\ref{bdn}) and (\ref{delta}),
respectively. Hence all conditions of Hayman's theorem hold and we
can apply it with $G_n=p(n), G(x)=g(x), r_n=e^{-h_n}$ and $\rho=1$
to find that
\begin{equation}\label{asypn}
p(n)\sim\frac{e^{nh_n}g(e^{-h_n})}{\sqrt{2\pi b(e^{-h_n})}}, \quad
n\to\infty,
\end{equation}
which completes the proof. $\rule{2mm}{2mm}$

{\it Remark.} To show that formula (\ref{asypn}) yields
(\ref{hardy}), one has to replace (\ref{dn}) and (\ref{bdn}) in
the right hand side of (\ref{asypn}).
 The asymptotic of $g(e^{-h_n})$ is
determined by a general lemma due to Meinardus [11] (see also [2;
Lemma 6.1]). Since $\zeta(0)=-1/2$ and
$\zeta^\prime(0)=-\frac{1}{2}\log{(2\pi)}$ (see [1; Chapter
23.2]), in the particular case of $g(e^{-h_n})$ this lemma implies
that
\begin{eqnarray}
& & g(e^{-h_n}) =\exp{(\zeta(2)h_n^{-1} -\zeta(0)\log{h_n}
+\zeta^\prime(0) +O(h_n^{c_1}))} \nonumber \\
& & =\exp{\left(\frac{\pi^2}{6h_n} +\frac{1}{2}\log{h_n}
-\frac{1}{2}\log{(2\pi)} +O(h_n^{c_1})\right)}, \quad n\to\infty,
\nonumber
\end{eqnarray}
 where $0<c_1<1$. The rest of the computation leading to (\ref{hardy}) is based on simple
 algebraic manipulations and cancellations.

\section{Proof of Theorem 1}
 \setcounter{equation}{0}

 We base our proof on the definition of Sampling Procedure 1 and
 eq. (\ref{procone}). We want to replace the expected value in
 its right-hand side by the ratio
 $\mathbb{E}(Z_{d,s}^{(n)})/\mathbb{E}(Z_n)$. So, we notice first
 that Erd\"{o}s and Lehner [4] proved that, in probability, the
 total number of parts $Z_n$ is asymptotic to $\mathbb{E}(Z_n)$ as
 $n\to\infty$. Hence, for any $\epsilon>0$, the probability of the event
 $$
 A_n=\left\{\lambda\in\Lambda(n): \mid\frac{Z_n}{\mathbb{E}(Z_n)}-1\mid>\epsilon\right\}
 $$
 tends to $0$ as $n\to\infty$. Further, we rewrite (\ref{procone})
 in the following way:
 \begin{equation}\label{indicators}
 \mathbb{P}(\mu_{n,1}\le d, \sigma_{n,1}\le s)
=\mathbb{E}\left(\frac{Z_{d,s}^{(n)}}{Z_n}I_{A_n^c}\right)
+\mathbb{E}\left(\frac{Z_{d,s}^{(n)}}{Z_n}I_{A_n}\right).
\end{equation}
For $\lambda\in A_n^c$ and $0<\epsilon<1$, we have
$(1-\epsilon)\mathbb{E}(Z_n)\le Z_n\le
(1+\epsilon)\mathbb{E}(Z_n)$ and therefore,
\begin{equation}\label{twoineq}
\frac{\mathbb{E}(Z_{d,s}^{(n)})}{(1+\epsilon)\mathbb{E}(Z_n)} \le
\mathbb{E}\left(\frac{Z_{d,s}^{(n)}}{Z_n} I_{A_n^c}\right) \le
\frac{\mathbb{E}(Z_{d,s}^{(n)})}{(1-\epsilon)\mathbb{E}(Z_n)}.
\end{equation}
Since $Z_{d,s}^{(n)}\le Z_n$, the second summand in
(\ref{indicators}) is not greater than $\mathbb{P}(A_n)$. Hence,
combining (\ref{indicators}) and (\ref{twoineq}), we obtain
$$
\mathbb{P}(\mu_{n,1}\le d, \sigma_{n,1}\le s)
=(1+O(\epsilon))\frac{\mathbb{E}(Z_{d,s}^{(n)})}{\mathbb{E}(Z_n)}
+\mathbb{P}(A_n).
$$
Letting $n\to\infty$ and then $\epsilon\to 0$ and replacing
$\mathbb{E}(Z_n)$ by the right-hand side of (\ref{ezn}), uniformly
for $d,s\ge 1$, we finally get
\begin{equation}\label{firststart}
\mathbb{P}(\mu_{n,1}\le d, \sigma_{n,1}\le s)
\sim\frac{2\pi\mathbb{E}(Z_{d,s}^{(n)})}{\sqrt{6n}\log{n}}
=\frac{2c\mathbb{E}(Z_{d,s}^{(n)})}{\sqrt{n}\log{n}},
\end{equation}
where $c$ is the constant from (\ref{c}).

Our proof continues with an application of Cauchy coefficient
formula to (\ref{idezds}). We use the
 circle $x=e^{-h_n+i\theta}, -\pi<\theta\le\pi$, as a contour of
 integration and the notation
 \begin{equation}\label{phids}
 \varphi_{d,s}(x)= \left(\sum_{1\le j\le s}\frac{x^j}{1-x^j}
 -(d+1)\sum_{1\le j\le s} \frac{x^{j(d+1)}}{1-x^{j(d+1)}}\right)
 \prod_{1\le j\le s} (1-x^{j(d+1)})
 \end{equation}
 to obtain
  $$
 p(n)\mathbb{E}(Z_{d,s}^{(n)}) =\frac{e^{nh_n}}{2\pi}\int_{-\pi}^\pi
 g(e^{-h_n+i\theta}) \varphi_{d,s}(e^{-h_n+i\theta}) e^{-i\theta n}d\theta.
 $$
 Then, we break up the range of integration as follows:
 \begin{equation} \label{sumint}
p(n)\mathbb{E}(Z_{d,s}^{(n)}) =J_1(d,s,n) +J_2(d,s,n),
 \end{equation}
 where
 \begin{equation}\label{jone}
 J_1(d,s,n) =\frac{e^{nh_n}}{2\pi}\int_{-\delta_n}^{\delta_n}
g(e^{-h_n+i\theta}) \varphi_{d,s}(e^{-h_n+i\theta}) d\theta,
\end{equation}
\begin{equation} \label{jtwo}
J_2(d,s,n)
=\frac{e^{nh_n}}{2\pi}\int_{\delta_n<\mid\theta\mid\le\pi}
g(e^{-h_n+i\theta}) \varphi_{d,s}(e^{-h_n+i\theta}) d\theta
\end{equation}
and $\delta_n$ is defined by (\ref{delta}).

In our next step we set
\begin{equation}\label{dands}
d=n^{u/2}, \quad s=n^{v/2}, \quad 0\le u,v\le 1
\end{equation}
and obtain estimates for the sums:
\begin{equation}\label{esone}
S_1=\sum_{1\le j\le s}\frac{e^{-jh_n}}{1-e^{-jh_n}},
\end{equation}
\begin{equation}\label{estwo}
S_2=\sum_{1\le j\le s}\frac{e^{-j(d+1)h_n}}{1-e^{-j(d+1)h_n}}.
\end{equation}
Here the sequence $\{h_n\}_{n\ge 1}$ is defined by (\ref{dn}).

Using the approximation of a Riemann sum by an integral,
(\ref{dn}), (\ref{dands}) and (\ref{c}), for $S_1$ we get
\begin{eqnarray}\label{esoneas}
& & S_1=\left(1+O\left(\frac{1}{n}\right)\right)
\sqrt{n}\sum_{1\le j\le n^{v/2}}
\frac{e^{-cj/\sqrt{n}}}{1-e^{-cj/\sqrt{n}}}\frac{1}{\sqrt{n}}
 \nonumber \\
 & & \sim\sqrt{n}\int_{1/\sqrt{n}}^{n^{\frac{v-1}{2}}}
 \frac{e^{-cz}}{1-e^{-cz}}dz =\frac{\sqrt{n}}{c} \int_{c/\sqrt{n}}^{n^{\frac{v-1}{2}}}
 \frac{e^{-z}}{1-e^{-z}}dz \nonumber \\
 & & =\frac{\sqrt{n}}{c}
 \log{\left(\frac{1-e^{-cn^{\frac{v-1}{2}}}}{1-e^{-c/\sqrt{n}}}\right)}
 =\frac{v\sqrt{n}}{2c}\log{n} +O(n^{v/2}).
\end{eqnarray}
In the same way one can show that
\begin{equation}\label{estwoas}
S_2= \left\{\begin{array}{ll} \frac{1-u}{2c} n^{\frac{1-u}{2}}\log{n}
+O(n^{\frac{1-u}{2}}) & \qquad  \mbox {if}\qquad v+u\ge 1, \\
 \frac{v}{2c} n^{\frac{1-u}{2}}\log{n} +O(n^{\frac{1-u}{2}}) &
\qquad \mbox {if}\qquad v+u<1.
 \end{array}\right.
 \end{equation}

We are now ready to find an estimate for the second integral in
(\ref{sumint}) (see (\ref{jtwo})). First, we have
\begin{eqnarray}\label{three}
& & \mid\prod_{1\le j\le s} (1-e^{-h_n j(d+1)+ij\theta(d+1)}\mid
\nonumber \\
& & \le \prod_{1\le j\le s} (1-e^{-h_n j(d+1)}) +\prod_{1\le j\le
s}
e^{-h_n j(d+1)}\mid 1-e^{ij\theta(d+1)}\mid \nonumber \\
& & \le 1+e^{-h_n(d+1)s}(1+\prod_{1\le j\le s}\mid
e^{ij\theta(d+1)}\mid) =1+2e^{-h_n(d+1)s} \le 3.
\end{eqnarray}
Hence, in terms of notations (\ref{phids}), (\ref{esone}) and
(\ref{estwo}), by (\ref{dands}), (\ref{esoneas}) and
(\ref{estwoas}),
$$
\mid\varphi_{d,s}(e^{-h_n+i\theta})\mid =O((S_1+(d+1)S_2)
=O(\sqrt{n}\log{n}), \quad -\pi\le\theta\le\pi.
$$
Replacing this estimate and applying inequality (\ref{decay}) to
the integrand of (\ref{jtwo}), we obtain
 $$
 \mid J_2(d,s,n)\mid =O(e^{n h_n}g(e^{-h_n})\sqrt{n}(\log{n})
 e^{-c_0 h_n^{-\epsilon_0}}).
 $$
 The required estimate now follows from (\ref{bdn}) and
 (\ref{asypn}) in the following way:
 \begin{eqnarray}\label{jtwoas}
 & & \mid J_2(d,s,n)\mid =O\left(\frac{e^{n h_n}g(-e^{-h_n})}
 {\sqrt{b(e^{-h_n}})}n^{1/2+3/4}(\log{n})e^{-c_0
 h_n^{-\epsilon_0}}\right) \nonumber \\
 & & =O(p(n)n^{5/4}(\log{n})e^{-c_0
 h_n^{-\epsilon_0}}) =O(p(n)n^{5/4}(\log{n})e^{-c_0^\prime
 n^{\epsilon_0/2}}) \nonumber \\
 & & =o(\sqrt{n}(\log{n})p(n)),
 \end{eqnarray}
where $c_0^\prime>0$.

The estimate for $J_1(d,s,n)$ follows from Hayman's "locality"
condition (\ref{locality}). First, we need to expand
$\varphi_{d,s}$ by Taylor formula. We have
\begin{eqnarray}\label{taylor}
& & \varphi_{d,s}(e^{-h_n+i\theta}) =\varphi_{d,s}(e^{-h_n})+
O\left(\mid\theta\mid\frac{d}{dx}\varphi_{d,s}(x)\mid_{x=e^{-h_n}}\right)
\nonumber \\
& & =\varphi_{d,s}(e^{-h_n})
+O\left(\delta_n\frac{d}{dx}\varphi_{d,s}(x)\mid_{x=e^{-h_n}}\right).
\end{eqnarray}
To find the asymptotic of $\varphi_{d,s}(e^{-h_n})$, in addition
to (\ref{esoneas}) and (\ref{estwoas}), we also need the limit of
$\prod_{1\le j\le s}(1-e^{-j(d+1)h_n})$ as $n\to\infty$, whenever
$d$ and $s$ satisfy (\ref{dands}) (see (\ref{phids})). Using
approximations by Riemann integrals as in the analysis of $S_1$
and $S_2$, it is easy to show that
$$
\prod_{1\le j\le s}(1-e^{-j(d+1)h_n}) =\exp{\left(\sum_{1\le j\le
s}\log{(1-e^{-j(d+1)h_n})}\right)} \to \left\{\begin{array}{ll} 1
& \qquad
\mbox {if}\qquad v+u\ge 1, \\
0 & \qquad \mbox {if} \qquad v+u<1.
 \end{array}\right.
 $$
Hence, from (\ref{dands})-(\ref{estwoas}) it follows that
\begin{equation}\label{taylorone}
\varphi_{d,s}(e^{-h_n}) \left\{\begin{array}{ll}
\sim\frac{u+v-1}{2c} \sqrt{n}\log{n}
 & \qquad  \mbox {if}\qquad v+u\ge 1, \\
 =o(\sqrt{n}) &
\qquad \mbox {if}\qquad v+u<1.
 \end{array}\right.
 \end{equation}
 The estimate of the error term in (\ref{taylor}) is tedious and
 follows the same line of reasoning. We have
 \begin{eqnarray}\label{error}
 & & \frac{d}{dx}\varphi_{d,s}(x)\mid_{x=e^{-h_n}} \nonumber \\
 & & =\left(\sum_{1\le j\le s}\frac{jx^{j-1}}{(1-x^j)^2}
 -(d+1)\sum_{1\le j\le
 s}\frac{jx^{j(d+1)}}{(1-x^{j(d+1)})^2}\right) \nonumber \\
 & & \times\left(\prod_{1\le j\le s} (1-x^{j(d+1)})\right)\mid_{x=e^{-h_n}}
 +\left(\prod_{1\le j\le s} (1-x^{j(d+1)})\right) \nonumber \\
 & & \times\exp{\left(-(d+1)\sum_{1\le j\le s}\frac{jx^{j(d+1)-1}}{1-x^{j(d+1)}}\right)}\mid_{x=e^{-h_n}}.
  \end{eqnarray}
It can be seen that the first two sums in the right-hand side of
(\ref{error}) are of order $O(n\log{n})$, while the first product
factor is estimated by (\ref{three}). Hence, the first summand in
(\ref{error}) is of order $O(n\log{n})$. For the sum in the
exponent of the second summand of the right-hand side of
(\ref{error}), one can show that there exists a constant $C>0$
such that
$$
(d+1)\sum_{1\le j\le
s}\frac{jx^{j(d+1)-1}}{1-x^{j(d+1)}}\mid_{x=e^{-h_n}}\ge
C\sqrt{n}\log{n}.
$$
Therefore the second summand in (\ref{error}) is
$O(e^{-C\sqrt{n}\log{n})})$. Hence
$$
\frac{d}{dx}\varphi_{d,s}(x)\mid_{x=e^{-h_n}} =O(n\log{n})
$$
and by (\ref{delta}) and (\ref{taylorone}), the expansion in
(\ref{taylor}) becomes
 \begin{equation}\label{phidsasym}
 \varphi_{d,s}(e^{-h_n+i\theta}) =\varphi_{d,s}(e^{-h_n})
+O\left(\frac{n^{1/3}\log{n}}{\Omega(n)}\right),
 \end{equation}
where $\Omega(n)\to\infty$ as $n\to\infty$ arbitrarily slowly.
Inserting this estimate and (\ref{locality}) into (\ref{jone}) and
applying the asymptotic for the partition function $p(n)$ from
(\ref{asypn}), we obtain
\begin{eqnarray}\label{joneest}
& & J_1(d,s,n) =\frac{e^{nh_n}g(e^{-h_n})}{2\pi}
\left(\int_{-\delta_n}^{\delta_n} e^{-\theta^2 b(e^{-d_n})/2}
(1+O(1/\Omega^3(n))d\theta\right) \nonumber \\
& & \times\left(\varphi_{d,s}(e^{-h_n})
+O\left(\frac{n^{1/3}\log{n}}{\Omega(n)}\right)\right) \nonumber \\
& & \sim\frac{e^{nh_n}g(e^{-h_n})}{\sqrt{b(e^{-h_n})}2\pi}
\left(\int_{-\delta_n\sqrt{b(e^{-h_n})}}^{\delta_n\sqrt{b(e^{-h_n})}}
e^{-y^2/2}dy\right) \left(\varphi_{d,s}(e^{-h_n})
+O\left(\frac{n^{1/3}\log{n}}{\Omega(n)}\right)\right) \nonumber
\\
& & \sim\frac{e^{nh_n}g(e^{-h_n})}{\sqrt{b(e^{-h_n})}2\pi}
\left(\int_{-\infty}^\infty e^{-y^2/2}dy\right)
\left(\varphi_{d,s}(e^{-h_n})
+O\left(\frac{n^{1/3}\log{n}}{\Omega(n)}\right)\right) \nonumber
\\
& & =\frac{e^{nh_n}g(e^{-h_n})}{\sqrt{2\pi
b(e^{-h_n})}}\left(\varphi_{d,s}(e^{-h_n})
+O\left(\frac{n^{1/3}\log{n}}{\Omega(n)}\right)\right)
\nonumber \\
& & \sim p(n)\left(\varphi_{d,s}(e^{-h_n})
+O\left(\frac{n^{1/3}\log{n}}{\Omega(n)}\right)\right),
\end{eqnarray}
where for the second asymptotic equivalence we have used
(\ref{bdn}) and (\ref{delta}) in order to get
$$
\delta_n\sqrt{b(e^{-h_n})}\sim
\frac{\pi^{5/6}\sqrt{2}}{6^{1/6}\Omega(n)} n^{1/12}\to\infty
$$
if $\Omega(n)\to\infty$ as $n\to\infty$ not too fast, so that
$\frac{n^{1/12}}{\Omega(n)}\to\infty$. It is now clear that
(\ref{sumint})-(\ref{jtwo}), (\ref{jtwoas}) and (\ref{joneest})
yield
$$
p(n)\mathbb{E}(Z_{d,s}^{(n)}) =p(n)\varphi_{d,s}(e^{-h_n})
+o(p(n)\sqrt{n}\log{n})
$$
and therefore
$$
\mathbb{E}(Z_{d,s}^{(n)}) =\varphi_{d,s}(e^{-h_n})
+o(\sqrt{n}\log{n}).
$$
 The result of Theorem 1 now follows from
(\ref{firststart}), (\ref{ezn}), (\ref{c}) and (\ref{taylorone}).

\section{Proof of Theorem 2}
\setcounter{equation}{0}

 We base our proof on (\ref{proctwo}),
Lemmas 2 and 3 and asymptotic equivalence (\ref{eyn}). To replace
the expectation in the right hand side of (\ref{proctwo}) by the
ratio $\mathbb{E}(Y_{m,s}^{(n)})/\mathbb{E}(Y_n)$, similarly to
what we did in the proof of Theorem 1, we shall study how unlikely
is the event
$$
B_n=\left\{\lambda\in\Lambda(n):
\mid\frac{cY_n(\lambda)}{\sqrt{n}}-1\mid>\epsilon\right\}, \quad
\epsilon>0,
$$
where $c$ is constant from (\ref{c}). Using Fristeft's method [6],
Corteel et al. [3] showed that
\begin{equation}\label{cortbound}
Pr(B_n)\le e^{-c_2\sqrt{n}}, \quad c_2=c_2(\epsilon)>0.
\end{equation}

{\it Remark.} Fristedt's approach [6] is based on the identity
\begin{equation}\label{probabgeom}
Pr(\alpha_j^{(n)}=m_j, j=1,...,n) =Pr\left(\gamma_j=m_j,
j=1,...,n\mid\sum_{j\ge 1} j\gamma_j=n\right),
\end{equation}
where $\{\gamma_j\}_{j\ge 1}$ is a sequence of independent
geometrically distributed random variables, whose distribution is
given by
$$
Pr(\gamma_j=k) =(1-q^j)q^{jk}, \quad k=0,1,...
$$
and $\{m_j\}_{j\ge 1}$ are non-negative integers. Eq.
(\ref{probabgeom}) holds for every fixed $q\in (0,1)$. It is
natural to take $q$ so that $Pr(\sum_{j\ge 1}j\gamma_j=n)$ is as
large as possible. Fristedt's almost optimal choice for $q$ is
$q=e^{-c/\sqrt{n}}$. Then, the bound in (\ref{cortbound}) is
easily obtained using this value of $q$.

Next, we represent the probability in (\ref{proctwo}) in the
following way
\begin{equation}\label{ind}
\mathbb{P}(\mu_{n,2}=m, \sigma_{n,2}\le s)
=\mathbb{E}\left(\frac{Y_{m,s}^{(n)}}{Y_n}I_{B_n^c}\right)
+\mathbb{E}\left(\frac{Y_{m,s}^{(n)}}{Y_n}I_{B_n}\right),
\end{equation}
where $I_{B_n}$ and $I_{B_n^c}$ denote the indicators of events
$B_n$ and $B_n^c$, respectively. Since, for any $\lambda\in
B_n^c$,
$$
\frac{c}{\sqrt{n}(1+\epsilon)}<\frac{1}{Y_n}<\frac{c}{\sqrt{n}(1-\epsilon)}
$$
if $0<\epsilon<1$, the first summand in (\ref{ind}) is estimated
by
\begin{eqnarray}\label{biensi}
& & \mathbb{E}\left(\frac{Y_{m,s}^{(n)}}{Y_n}I_{B_n^c}\right)
=\frac{c}{\sqrt{n}}(1+O(\epsilon))\mathbb{E}(Y_{m,s}^{(n)}I_{B_n^c})
\nonumber \\
& & =\frac{c}{\sqrt{n}}(1+O(\epsilon))(\mathbb{E}(Y_{m,s}^{(n)})
-\mathbb{E}(Y_{m,s}^{(n)}I_{B_n})).
\end{eqnarray}
Clearly, with probability $1$, $Y_{m,s}^{(n)}\le n$. Hence, using
(\ref{cortbound}), we obtain
$$
\mathbb{E}(Y_{m,s}^{(n)}I_{B_n})) =O(nPr(B_n))
=O(ne^{-c_2\sqrt{n}})
$$
and (\ref{biensi}) becomes
$$
\mathbb{E}\left(\frac{Y_{m,s}^{(n)}}{Y_n}I_{B_n^c}\right)
=\frac{c}{\sqrt{n}}(1+O(\epsilon))\mathbb{E}(Y_{m,s}^{(n)})
+O(ne^{-c_2\sqrt{n}}).
$$
The second term in the right hand side of (\ref{ind}) is easily
estimated using (\ref{cortbound}) since it is not greater than
$Pr(B_n)$. Consequently,
\begin{equation}\label{twotwo}
\mathbb{P}(\mu_{n,2}=m,\sigma_{n,2}\le s)
=\frac{c}{\sqrt{n}}(1+O(\epsilon))\mathbb{E}(Y_{m,s}^{(n)})
+O(ne^{-c_2\sqrt{n}})
\end{equation}
uniformly for any fixed integer $m\ge 1$ and real $s\ge 1$. Hence,
our next task is to obtain an estimate for
$\mathbb{E}(Y_{m,s}^{(n)})$, as $n\to\infty$, whenever
$s=t\sqrt{n}/c$, $m\ge 1$ is fixed integer and $t\in (0,\infty)$
is also fixed. Combining results of (\ref{eyms}) and (\ref{ratio})
of Lemmas 2 and 3, respectively, and approximating the sum by the
corresponding Riemann integral, we get
\begin{eqnarray}
& & \mathbb{E}(Y_{m,s}^{(n)}) =(1+O(1/\sqrt{n}))\nonumber \\
& & \times\sum_{1\le k\le c^{-1}\sqrt{n}t} (\exp{(-cmk/\sqrt{n})}
-\exp{(-c(m+1)k/\sqrt{n})})
\nonumber \\
& & \sim \sqrt{n}\int_0^{c^{-1}t} (e^{-cmy}-e^{-c(m+1)y})dy.
\nonumber
\end{eqnarray}
 Replacing this expression into (\ref{twotwo}) and
letting first $n\to\infty$ and then $\epsilon\to 0$, we obtain
\begin{eqnarray}
& & \mathbb{P}\left(\mu_{n,2}=m,\frac{c\sigma_{n,2}}{\sqrt{n}}\le
t\right) \to c\int_0^{c^{-1}t}
(e^{-cmy}-e^{-c(m+1)y})dy \nonumber \\
& & =\int_0^t e^{-my}(1-e^{-y})dy, \nonumber
\end{eqnarray}
which completes the proof of Theorem 2.

\section{Proof of Theorem 3}
\setcounter{equation}{0}

The proof will be based on an asymptotic analysis of formula
(\ref{procthree}), setting there $s=c^{-1}\sqrt{n}t$ as
$n\to\infty$ (see again (\ref{c})) and assuming that $m$ is fixed
positive integer. First, we let $\Lambda_k(n)$ to denote the set
of partitions of $n$ with no part equal to $k$. Also, let
$P_k(n)=\mid\Lambda_k(n)\mid$. In [3; p. 189] Corteel et al. give
a combinatorial proof of the following identity:
$$
Pr(\alpha_j^{(n)}=m) =\frac{P_j(n-mj)}{p(n)}
=\frac{p(n-jm)-p(n-j(m+1))}{p(n)}.
$$
Replacing this expression into the right hand side of
(\ref{procthree}) and applying (\ref{ratio}), as in the proof of
Theorem 2, we obtain
\begin{eqnarray}
& & \mathbb{P}\left(\mu_{n,3}=m,\frac{c\sigma_{n,3}}{\sqrt{n}}\le
t\right) \sim m\sum_{1\le j\le c^{-1}\sqrt{n}t}\frac{j}{\sqrt{n}}
(e^{-cmj/\sqrt{n}}-e^{c(m+1)j/\sqrt{n}})\frac{1}{\sqrt{n}}
\nonumber \\
& & \to m\int_0^{c^{-1}t} y(e^{-cmy}-e^{-c(m+1)y})dy
=\frac{m}{c^2}\int_0^t y(1-e^{-y})e^{-my}dy. \nonumber
\end{eqnarray}
This completes the proof of Theorem 3.

\end{document}